\documentclass[a4paper,11pt]{article}
\usepackage[utf8]{inputenc}
\usepackage{appendix}
\usepackage[margin=3cm]{geometry}
\usepackage{amsfonts,amsmath,amssymb,amscd,graphicx,mathabx}
\usepackage[colorlinks]{hyperref}
\usepackage{enumitem}

\usepackage{blkarray}

\newtheorem{theorem}{Theorem}[section]
\newtheorem{corollary}[theorem]{Corollary}
\newtheorem{lemma}[theorem]{Lemma}

\newtheorem{proposition}[theorem]{Proposition}

\newtheorem{definition}[theorem]{Definition}
\newtheorem{remark}[theorem]{Remark}
\newtheorem{example}[theorem]{Example}
\numberwithin{equation}{section}

\newenvironment{preuve}[1][]
{\vskip 2mm  \emph{\bf Proof#1. }}{$\Box$ \vskip 2mm}

\newcommand{\bpr}{\begin{preuve}}
\newcommand{\epr}{\end{preuve}}

\newcommand{\beq}{\begin{eqnarray}}
\newcommand{\eeq}{\end{eqnarray}}
\newcommand{\beqe}{\begin{eqnarray*}}
\newcommand{\eeqe}{\end{eqnarray*}}

\newcommand{\Nn}{\mathbb{N}}
\newcommand{\R}{\mathbb{R}}

\newcommand{\C}{\mathbb{C}}

\let\epsilon=\varepsilon

\newcommand{\vol}{\ensuremath \mathrm{vol}}
\newcommand{\Tr}{\ensuremath \mathrm{Tr}}

\newcommand{\crit}{\ensuremath \mathrm{Crit}}
\newcommand{\Ind}{\ensuremath \mathrm{Ind \, }}
\newcommand{\grass}{\ensuremath \mathrm{Grass}}
\newcommand{\cov}{\ensuremath \mathrm{Cov}}
\newcommand{\sym}{\ensuremath \mathrm{Sym}}
\newcommand{\End}{\ensuremath \mathrm{End}}
\newcommand{\id}{\ensuremath \mathrm{Id}}
\begin{document}

\title{\bf Expected local topology of \\random complex submanifolds}

\author{\sc Damien Gayet}
\maketitle

\begin{abstract}
Let $n\geq 2$ and $r\in \{1, \cdots, n-1\}$ be integers, $M$ be a compact smooth K\"ahler manifold of complex dimension $n$, $E$ be a holomorphic vector bundle with complex rank $r$ and equipped with an hermitian metric $h_E$, and $L$ be an ample holomorphic line bundle over $M$ equipped with a metric $h$ with positive curvature form. For any $d\in \mathbb N$ large enough, we equip the space of holomorphic sections $H^0(M,E\otimes L^d)$ with the natural Gaussian measure associated to $h_E$ , $h$ and its curvature form. Let $U\subset M$ be an open subset with smooth boundary. We prove that the average of the $(n-r)$-th Betti number of the vanishing locus in $U$ of a random section 
$s$  of $H^0(M,E\otimes L^d)$ is asymptotic to ${n-1 \choose r-1} d^n\int_U c_1(L)^n$ for large $d$. On the other hand, the average of the other Betti numbers are $o(d^n)$. The first asymptotic recovers the classical deterministic global algebraic computation. Moreover, such a discrepancy in the order of growth of these averages is new and constrasts  with all known other smooth Gaussian models, in particular the real algebraic one. We prove a similar result for the affine complex Bargmann-Fock model.
\end{abstract}

\textsc{Mathematics subject classification  2010}: 32L05 (Holomorphic bundles and generalizations) ; 60G15 (Gaussian processes)\\

\tableofcontents

\section{Introduction}
The goal of this article is to understand the statistics of the local topology 
of random complex submanifolds, for projective manifolds and the affine complex space.
\paragraph{Projective manifolds.}
Let $n$ be a positive integer, $M$ be a compact smooth complex manifold of complex dimension $n$, and $L$ be an ample holomorphic line bundle over $M$. Let $h$ be a Hermitian metric on $L$ with positive curvature form $c_1(L) =\omega$, that is locally 
\beq\label{omega} \omega = \frac{1}{2i\pi}\partial \bar\partial \log \|s\|_h^2,
\eeq
where $s$ is any local non-vanishing section of $L$.
Then, $(M,\omega)$ becomes a K\"ahler manifold, and by the Kodaira theorem, it can be embedded in a projective space.
For any large enough degree $d\geq 1$, and any generic holomorphic section $s\in H^0(M, L^d)$, denote by $Z_s\subset M$ the smooth vanishing locus of $s$. The famous hyperplane Lefschetz theorem asserts, in particular, that~\cite{griffiths} 
$$\forall 0\leq i\leq n-2, \ b_i(Z_s)=b_i(M).$$ 
For instance, if $M=\C P^n$, then for $i\leq n-2$, $b_i(Z_s)=0$ if $i$ is odd and $b_i(Z_s)=1$ if $i$ is even. 
On the other hand, 
\beq\label{asymp}
\frac{1}{d^n} b_{n-1}(Z_s)\underset{d\to \infty}{\to} \int_M \omega^n.
\eeq	
Of course, there is no local (deterministic) version of the  Lefschetz theorem. Indeed, if $U$ is an open subset of $M$,  the intersection of $U$ with $Z_s$ can be empty or can have a topologicial complexity bigger  than the one of $Z_s$. In particular for $n\geq 2$, $Z_s$ is connected but its intersection with $U$ can be disconnected. 
There is even no bound for the number of components of it, since we can twist $U$ for that.  However, for a fixed $U$ defined by algebraic inequalities, the following bound exists:
\begin{theorem}(\cite[Theorem 3]{milnor}) Let $U\subset \C P^n\setminus\{Z_0=0\}$ be an open subset defined by real algebraic inequalities. Then, there exists a constant $C_U$ depending only on the number and the degree of the defining polynomials of $U$, such that for any generic $r$-uple of homogeneous complex polynomials $s=(s_1, \cdots, s_r)\in (\C^{hom}_d)^r$ of degree $d$, 
\beq\label{milnor}
\sum_{i=0}^{2n-2r} b_i(Z_s\cap U)\leq C_{U}d^{2n}.
\eeq
\end{theorem}

Now, if the section $s$ is taken at random, one could hope that for fixed $U$, not necessarily defined by polynomials, the average topology of $Z_s\cap U$ reflects in some way the Lefschetz theorem and with further hope, the asymptotic~(\ref{asymp}) as well. In this paper, we prove that these two intuitions are true, in the following more general classical setting.
In addition to $(L,h)$, let $(E,h_E)$ be a holomorphic vector bundle of rank $r$ and equipped with a Hermitian metric $h_E$. Since $L$ is ample, for $d$ large enough, the space of holomorhic sections $H^0(M,E\otimes L^d)$ is non-trivial. 
Then, a natural scalar product associated to this setting is the following:
\beq\label{herm}
\forall (s,t)\in (H^0(M,E\otimes L^d))^2, \ 
\langle s,t\rangle = \int_M h_{E}\otimes h_{L^d}(s,t) \frac{\omega^n}{n!},
\eeq
where $h_{L^d}$ is the metric over $L^d$ induced by $h^d$. 
A natural probability measure $\mu_d$ over this space is the Gaussian one associated to this Hermitian product. 
In other terms, for any Borelian $A\subset H^0(M,E\otimes L^d)$, 
\beq\label{mesure}
\mu_d( A) = \int_{A}e^{-\frac12 \|s\|^2} \frac{ds}{(2\pi)^{N_d}},
\eeq
where $\|\cdot\|$ denotes the norm associated to the Hermitian product~(\ref{herm}), $N_d$ the (complex) dimension of $H^0(M,E\otimes L^d)$ and $ds$ the Lebesgue measure. 
Notice that if $(S_i)_{i\in\{1, \cdots, N_d\}}$ is an orthonormal basis for this scalar product, 
	then $s=\sum_{i=1}^{N_d} a_i S_i$ is random for $\mu_d$
	when the coefficients $a_i\in \C$ are i.i.d standard complex Gaussians, that is $\Re a_i $ and $\Im a_i$ are independent standard Gaussians.
	\begin{example}\label{exa}
		For $M=\C P^n$, $E=\C^r$ equipped with its standard Hermitian metric, $L=\mathcal O(1)$ equipped with the Fubini-Studi metric, then $s$ consists in $r$ independent copies of random polynomials 
		$$ \forall 1\leq i\leq r, \ 
		s_i([Z])=\sum_{i_0+\cdots +i_n =d}a_{i_0\cdots i_n}\sqrt{\frac{(n+d)!}{n!i_0!\cdots i_n!}} Z_0^{i_0}\cdots Z_n^{i_n},$$
		where the $(a_i)_i$ are independent standard complex Gaussian variables. 
	\end{example}
Our main result is the following:

	\begin{theorem}\label{main} 
		Let $n\geq 2$ and $1\leq r\leq n-1$  be integers, $M$ be a compact smooth K\"ahler manifold and $(L,h)$ be an ample complex line bundle over $M$, with positive curvature form $\omega$, $(E,h_E)$ be a holomorphic rank $r$ vector bundle and
let		$U\subset M$ be a 
		open subset with smooth boundary. Then
		\begin{eqnarray*}
			\forall i\in \{0,\cdots, 2n-2r\}\setminus\{n-r\}, \ \frac{1}{d^n}\mathbb E b_i(Z_s\cap U) &\underset{d\to \infty}{\to}& 0\\
			\frac{1}{d^n}\mathbb E b_{n-r} (Z_s\cap U)&\underset{d\to \infty}{\to} & {n-1 \choose r-1}\int_U \omega^n.
		\end{eqnarray*}	
	Here the probability measure is the Gaussian one given by~(\ref{mesure}). These asymptotics hold when $U=M$ as well.
	\end{theorem}
Of course, when $U=M$, the topological type of $Z_s$ does not depend on the random section $s$. 
Markov's inequality implies the following corollary.
\begin{corollary}\label{coro0} Under the hypotheses of Theorem~\ref{main}, for any $\epsilon>0$, 
	$$\limsup_{d\to +\infty}\mu_d\left\lbrace s\in H^0(M,E\otimes L^d) \ | \ 
	b_{n-r} (Z_s\cap U)\geq \frac{d^n}\epsilon {n-1 \choose r-1}\int_U \omega^n\right\rbrace\leq \epsilon,$$
	where $\mu_d$ is defined by~(\ref{mesure}).
	\end{corollary}
Note that the Gaussian measure can be replaced by the round metric on the sphere $\mathbb S H^0(M,E\otimes L^d)$, where the metric is defined by~(\ref{herm}). Hence, this corollary can
be seen as a deterministic result about the volume of certain
subsets of topological interest in this sphere. 
\begin{example} Under the standard setting of Example~\ref{exa}, 
	$\int_{\C P^n} \omega_{FS}^n = 1$, 
	so that 
$$	\frac{1}{d^n}\mathbb E b_{n-r} (Z_s\cap U)\underset{d\to \infty}{\to} {n-1 \choose r-1}\frac{\vol (U)}{\vol(\C P^n)} .$$
	\end{example}
\begin{remark}
	\begin{enumerate}
		\item Theorem~\ref{main} provides the first explicit asymptotic for one mean Betti numbers of the nodal set of a smooth Gaussian field. Former explicit asymptotics were proven~\cite{gayet21} in a real context for high level random sets (in particular, not the zero one). For $r=1$, it is striking that the asymptotic average local behaviour reflects exactly the global asymptotic estimate given by~(\ref{asymp}).
\item For $U=M$, Theorem~\ref{main} has a deterministic corollary. Indeed, any divisor of given degree is diffeomorphic to the another of the same degree. Hence, for $r=1$ the second assertion of Theorem~\ref{main} is equivalent to the asymptotic~(\ref{asymp}). For higher codimensions $r$, \cite[Corollary 3.5.2]{gayet2015expected} shows that
\beq\label{chi}
\frac{1}{d^n}\chi(Z_s)\underset{d\to\infty}{\to} {n-1\choose r-1} \int_M \omega^n.
\eeq
Again, (\ref{chi}) is implied by Theorem~\ref{main}. 
\item 
We emphasize that these qualitatively different asymptotics are new. In particular constrasts  with the real situation~\cite[Corollary 1.2.2]{gayet2015expected} and all known others smooth Gaussian models like~\cite{GWasian} (see also~\cite{wigman2021expected}). In these latter cases, all Betti numbers grow like $L^n$, where $1/L$ is the natural scale of the model, $1/\sqrt d$ in this one. This is especially true for the number of connected components, see~\cite{nazarov2}. Here, the scale is $d^{-\frac12}$, however only the $(n-r)$-th Betti number grows like $d^n$.
\item In~\cite{gayet2019systoles}, it was proved that for any compact smooth real hypersurface $\Sigma$ of $\R^n$, for any open subset $U\subset M$,  with uniform probability, a uniform proportion of the $(n-1)$-homology in $Z_s\cap U$ can be represented by Lagrangians submanifolds diffeomorphic to $\mathcal L$.
	\item In \cite[Theorem 5 (2)]{ancona2020exponential} (see also~\cite[Theorem A]{diatta2021low}), it is shown that as far as (local) topology of $Z_s\cap \R P^n$ is only concerned, a random \emph{real} polynomial $s$ of degree $d$ can be replaced, with high probability, by a polynomial of degree slightly greater than $\sqrt d$. In fact, this statement holds for complex polynomials on a ball in the complementary of a complex hypersurface as well. Using Milnor's bound~(\ref{milnor}), this replacement allows to get a similar estimate as Corollary~\ref{coro0} when $U$ is be defined algebraically. The decay is almost exponential in this case.
\item  In~\cite[Proposition 6]{auroux}, 
the author proved that (deterministic) Donaldson hypersurfaces, which are zeros of sections with vanish transversally with a controlled derivative, satisfy such local topology estimate for the $(n-r)$-th Betti number. Theorem~\ref{main} shows a further evidence that Donaldson hypersurfaces have common features with random ones. For instance, the current of integration over $Z_s$ fills out uniformly $M$ for large degrees $d$ in both contexts, see~\cite{donaldson1996symplectic} and~\cite{shiffman}.
\end{enumerate}
\end{remark}

\paragraph{The complex Bargmann-Fock field.}
Finally, we prove an affine version in the universal limit for holomorphic sections, namely the complex Bargmann-Fock field.
The	Bargmann-Fock field is defined
by 
\beq\label{bbff}
\forall z\in \C^n, \ f(z)=\sum_{(i_1,\cdots, i_n)\in \Nn^n}a_{i_0, \cdots, i_n}
{\sqrt{\frac{\pi^{i_1+\cdots +i_n}}{i_1!\cdots i_n!}}}{z_1^{i_1}\cdots z_n^{i_n}}e^{-\frac12 \pi\|z\|^2},
\eeq
where the $a_I$'s are independent normal complex Gaussian random variables. The strange presence of $\pi$ will be explained below.

\begin{theorem}\label{BF}Let $n\geq 2$ and $1\leq r\leq n-1 $ be integers, $ f : \C^n\to \C^r$ be $r$ independent copies of the Bargmann-Fock field, and $U\subset \C^n$ be an open subset with compact smooth boundary. Then,
	\beqe 			\forall i\in \{0,\cdots, 2n-2r\}\setminus\{n-r\},\ 
	\ \frac{1}{R^{2n}}\mathbb E b_i( Z_f\cap RU)&\underset{R\to +\infty}\to& 0 \\
	\frac{1}{R^{2n}}	\mathbb E b_{n-r} ( Z_f\cap RU)&\underset{R\to +\infty}\to & n!{n-1 \choose r-1}\vol(U).
	\eeqe
	The volume is the standard one.
\end{theorem}
\begin{remark}
	\begin{enumerate}
\item Again, compared to the other known results, the order of magnitude of the mean number of connected components is not the natural one, that is $R^{2n}$, see~\cite{nazarov2} for instance. 
\item Theorem~\ref{BF} (and Theorem~\ref{main}) was guessed by the author for the following geometric reasons, which we present for $n=2$ and $r=1$: because of the maximum principle, if a complex curve in $\C^2$ locally touches a real hyperplane $H$, being (locally) on one side of $H$, then $C$ is affine and $C\subset H$. Now, if $p : U\subset \C^2\to \R$ is Morse, for any $R>0$, let $p_R= p(\cdot/R)$. Then  for large $R>0$, the level sets of $p_R$ are locally closer and closer to be planar, so that there should be less and less random cuves touching them from the interior, that is there are less and less critical points of $p_{|Z_f}$ of index 0, compared to critical points of index $1$. Morse theory should then imply the result.
\end{enumerate}
\end{remark}

\paragraph{Related results.}

The study of the statistics of the Betti numbers, or even the diffeomorphism type,  of a random smooth submanifold (of positive dimension) is now a well-developped subdomain of random geometry, with current links to percolation. We refer to~\cite{gayet21} for a historical account of this topic. The results were proven mainly in the real algebraic and Riemannian semiclassical settings. Both models share a common feature: the Betti numbers grow (with the parameter, degree or eigenvalue) like the inverse of the scale to a power equal to the dimension of the ambient manifold.   In both cases, the covariance of the model is the spectral kernel, for which estimates exist.

The local study of the geometry of random \emph{complex} submanifolds of positive dimension began with~\cite{shiffman}, under the hypotheses of Theorem~\ref{main}, with $r=1$. It was proven that the average current of integration over $Z_s$ tends to the curvature form of the line bundle, when $d$ grows to infinity. Since the topology of the complex hypersurfaces depend only on the degree, a crucial difference with the real setting, the topology of random complex hypersurfaces seemed less promising. Our paper~\cite{gayet2019systoles} showed that local random (symplectic) topology is interesting as well, and even can provide new deterministic results. 

A lot of results about critical points of random sections has been done. In this complex algebraic context, it seems to begin  with~\cite{douglas2004critical}. In~\cite{gayet2014total}, the restriction of a Lefschetz pencil to the complex random hypersurface was used in order to get topological estimates through Morse theory, which is the spirit of the present paper.
We refer to~\cite[\S 1.3]{gayet2015expected} for further references. 
The following result is close to the present work:
	\begin{theorem}\label{gw}{(\cite[Theorem 1.3]{GWdet} for $r=1$, \cite[Theorem 3.5.1]{gayet2015expected}) for any $r$) } 
	Under the hypotheses of Theorem~\ref{main}, let $p : M \dashrightarrow  \C P^{1}$
	be a Lefschetz pencil. 	Then,
	\begin{eqnarray*}
			\frac{1}{d^n}\mathbb E \#\left(U\cap \crit  (p_{|Z_s})\right) &\underset{d\to \infty}{\to} & {n-1 \choose r-1}\int_U \omega^n.
	\end{eqnarray*}	
\end{theorem}
This result holds in particular for any local holomorphic map.  A similar real version of Theorem~\ref{gw} was proven as well.
In the real setting, the authors used the weak Morse inequalities in order to get an upper bound for the average Betti numbers of $Z_s$. Lower bounds of the same order of magnitude (in the degree) where estimated by the barrier method. 

As in~\cite{gayet21}, in the present paper we use 
the \emph{strong} Morse inequalities, and moreover we use this theory on manifolds with boundary, which implies to take in account the critical points of the restriction of the function to the boundary. Joint with the weak ones, strong Morse inequalities allow us to get the proper estimate of the mean middle Betti number given by Theorem~\ref{main}. 
On the contrary to the real setting, in our complex setting strong Morse inequalities help, because the complex Hessian of a holomorphic function has a symmetric signature, which implies that all mean critical points of $p_{|Z_s}$ have the wrong order of magnitude, except when the index is the middle one, that is $n-r$, see Theorem~\ref{projective}.  

The method to prove Theorem~\ref{projective} is different than the one used for Theorem~\ref{gw}, but
both provide, on the one hand, a Kac-Rice formula (both based, at the end, on the coarea formula), and on the other hand, an estimate of it
when the degree goes to infinity. 
In~\cite{GWdet} and~\cite{gayet2015expected}, the authors used explicit peak sections to compute the average, and the aforementionned parts were mixed.
In this paper we wanted to clearly separate the two parts of the proof : one part which is based on a general Kac-Rice formula as Corollary~\ref{coro1}, and one part which depends on the particular model, real, holomorphic or mixed (on the boundary of the open set $U$). This can be done because the second part only needs informations about the covariance function. For projective manifolds, this is the Bergman kernel, see section~\ref{projo}. The peak sections is a way to recover the needed informations, see~\cite{tian1990set}. In~\cite{shiffman}, the Szeg\"o kernel was used, based on Zelditch's semiclassical way~\cite{zelditch1998szego} of proving Tian's theorem. For the Riemannian setting like in~\cite{GWasian}, the covariance is the spectral kernel and H\"ormander estimates can be used.

\paragraph{Holomorphic percolation.}
Theorem~\ref{BF} raises a natural question related to percolation theory: is there a Russo-Seymour-Welsh phenomenon for the complex Bargmann-Fock field ? In its simplest non-trivial form, this question is the following: \\

\noindent
{\it Let $B, B'\subset \mathbb S^3\subset \C^2$ two disjoint closed smooth $3-$balls lying in  the unit sphere, and let $f$ be the complex Bargmann-Fock field over $\C^2$ see~(\ref{bbff}). Is it true that}
$$\liminf_{R\to+\infty} \mathbb P\left(\exists \text{ a connected component of } \{f=0\}\cap R\mathbb B^4  \text{ joining } RB \text{ to } RB'\right)>0 \ ?$$
The analog for the real Bargmann-Fock over $\R^2$ is true, see~\cite{beffara2017percolation}. We emphasize that the holomorphic situations constrasts in many ways with the real setting. Firstly, there is no bounded component of $\{f=0\}$ in the complex case and with probability one, there is a unique component of $Z_f$. Secondly, none of the classical tools in percolation theory does hold in this holomorphic context, in particular duality and FKG property. Besides, the isotropy of the field and the absence of bounded components imply that with uniform positive probability there exists a component of $\{f=0\}\cap R\mathbb B$ from $\mathbb B$ to $RB$. In the real setting, this probability tends to 0. Finally, note that a similar question can be asked for complex algebraic submanifolds:

{\it Let $U\subset \C P^2$ be a smooth ball in the projective plane, $B, B'\subset \partial U$ two disjoint closed smooth $3-$balls lying in  the boundary of $U$, and let $s\in H^0(\C P^2,  \mathcal O(d))$ be a random polynomial of degree $d$. Is it true that }
	$$\liminf_{d\to+\infty} \mathbb P\left(\exists \text{ a connected component of } \{s=0\}\cap U  \text{ joining } B \text{ to } B'\right)>0 \ ?$$
	The real analog has been proven in~\cite{beliaev2021russo}.

\paragraph{Ideas of the proof of Theorem~\ref{main}.} Let $U\subset M$ be an open set with compact smooth boundary and $p : \bar U\to \R$ be a smooth Morse function in the sense of Definition~\ref{mobo}, that is $p$ Morse on $U$, its restriction to $\partial U$ is Morse and $p$ has no critical point on $\partial U$. Let $Z$ be a complex smooth submanifold of $U$ with boundary in $\partial U$, such that $p_{|Z}$ is Morse in the latter sense. Then, by Morse theory for manifolds with boundary, for any $0\leq i\leq \dim_\R Z$, the $i$-th Betti number of $Z$ is less or equal to the number of critical points of $p_{|Z}$ and $p_{|\partial Z}$ of index $i$, see Theorem~\ref{morse}. Besides, from the strong Morse inequalities, we can estimate the $i-$th Betti number of $Z$ if the critical points of index different than $i$ are far smaller. 

We apply this to $Z_s$ the zero set of a random holomorphic section $s$ of degree $d$. Note that the natural scale for the natural measure is $1/\sqrt d$. Hence, in every ball of this radius, the geometry of $Z_s$ should be independent of $d$. This implies in particular that on a manifold of real dimension $m$, the average of geometric or analytic observables like the number of critical points of $p_{|Z_z}$ should grow like $d^{\frac{m}2}$. 

The general Kac-Rice formula given by Corollary~\ref{coro1} applied to our projective situation allows us to estimate the number of critical points in the interior of $U$, see Theorem~\ref{BF0}, and a factor $d^n$ emerges, as guessed by the previous heuristic arguments. Here we use the fact that the covariance function of $s$ is the Bergman kernel and that this kernel has a universal rescaled limit, see Theorem~\ref{Dai}. Now, the integral in the Kac-Rice formula involves the determinant of a random matrix provided by a perturbation of the Hessian of $s$ (restricted to the tangent space of $Z_s$), where the perturbation decreases with $d$. At the limit, the matrix is non zero only for middle index, since the Hessian has complex symmetries. On the contrary, the mean number of critical points of middle index has a precise non-trivial asymptotic, see Theorem~\ref{projective}. 

We need also to control the number of critical points of the restriction of $p$ to the boundary of $U$ and of $Z_s$. We use the Kac-Rice formula in this mixed case as well, see Proposition~\ref{b0}. As guessed, a factor $d^{n-\frac12}$ emerges. Both estimates and Morse theory finish up the proof of Theorem~\ref{main}.

\paragraph{Structure of the article. }
In section~\ref{prel}, we prove various deterministic lemmas in order to prepare the main Kac-Rice formula computing the mean of critical points. This formula is established in section~\ref{kr}. In section~\ref{app}, we apply this formula in order to prove Theorem~\ref{main} and Theorem~\ref{BF}.

\paragraph{Aknowledgements.} The author thanks Thomas Letendre for his expertise of the book~\cite{ma2007holomorphic} and Michele Ancona for discussions about the subject of the article and his paper~\cite{ancona2020exponential}.

\section{Deterministic geometric preliminaries}\label{prel}

The general setting of this paper is a real manifold $M$ of dimension $n$ and a real vector bundle $E$ over $M$. Let also  $p: M\to \R$ be a Morse function. For any generic smooth section of $E$, we will look at the critical points of $p_{|\{s=0\}}$, which are the points $x\in M$ where 
the tangent space of the vanishing locus $Z_s$ of 
$s$ lies in $\ker dp(x)$. For this reason, we must understand the geometry of $\ker \nabla s$ as an element of the Grassmannian bundle $\grass(n-r,TM)$ or $\grass(n-r,\ker dp)$, where $\nabla$ denotes any covariant connection on $E$. In this section, we provide various simple lemmas which will be used in the main results. To make the computations easier, $M$ and $E$ will be endowed with metrics. 

\subsection{Kernel and Grassmannians}\label{kg}

The following Lemma is classical:
\begin{lemma}\label{gra}  
Let $n$ be an integer and $(V,g)$ be a finite dimensional real vector space of dimension $n$ equipped with a scalar product $g$.
For any integer $0\leq m\leq n$, the Grassmannian $\grass(m,V)$ of $m-$planes of $V$ is a smooth manifold and for any $K\in \grass(m,V)$, 
$T_K\grass(m,V)$ is canonically (with respect to $g$ and $K$) identified with $\mathcal L(K, K^\perp)$.
In particular, $T_K \grass(m,V)$ inherits the natural metric on $\mathcal L(K,K^\perp)$ induced by $g$. If $V$ is complex and $g$ is Hermitian, then the same holds, replacing the real Grassmaniann by the complex Grassmanian $\grass_\C$, and $ \mathcal L(K,K^\perp)$ by the complex linear maps $\mathcal L^\C (K,K^\perp)$. 
\end{lemma}
\begin{remark}Note that for $f\in \mathcal L(K,K^\perp)$, if $A$ denotes the matrix of $f$ in any $g$-orthonormal basis of $K$ and $K^\perp$, then the squared norm of $f$ induced by $g$ equals $\Tr (AA^*).$
	\end{remark}
\begin{lemma}\label{ram}Let $1\leq r\leq n$, $(V,g)$ be an Euclidean vector space, and $E$ be a real vector space, of respective dimensions $n$ and $r$. Let $\alpha_0 \in \mathcal L_{onto}(V,E)$ and $K=\ker \alpha_0\in \grass(n-r,V)$. Then, there exists a neighborhood $U\subset \mathcal L_{onto}(V,E)$ of $\alpha_0$ and a smooth 	map : $\varphi : U\to \mathcal L(K,K^\perp)$ such that 
	$\varphi(\alpha_0)=0$ and 
	$$\forall \alpha\in  U, \ \ker \alpha = \left( \id_{|K} +\varphi(\alpha)\right)(K).$$
Moreover, for any 
$  \beta\in \mathcal L(E,V), \  d\varphi(\alpha_0)(\beta)=- (\alpha_{0|K^\perp})^{(-1)}\beta_{|K}.$
The same holds in the complex-Hermitian setting. 
	\end{lemma} 
\bpr Let  \beqe F : \mathcal L(V,E)\times \mathcal L(K,K^\perp)&\to & \mathcal L(K,E)\\
(\alpha, f)& \mapsto & \alpha \circ (\ \id_{|K}+f).
\eeqe
Then, $F$ is smooth and $F(\alpha_0, 0)=0$.
The partial differential in $f$ at $(\alpha_0,0)$ writes
$$\forall g\in \mathcal L(K,K^\perp), \  d_f F(\alpha_0, 0)(g) = \alpha_0\circ g=\alpha_{0|K^\perp}\circ g\in \mathcal L(K,E).$$
This partial differential is an isomorphism because
$ \alpha_0$ is onto, so that $\alpha_{0|K^\perp} \in \mathcal L(K^\perp, E)$ is an isomorphism. Note that 
the partial differential in $\alpha$ satisfies
$$\beta\in \mathcal L(V,E),\   d_\alpha F (\alpha_0, 0) (\beta) = 
\beta_{|K}.$$
 Hence, by the implicit function theorem, there are two open neighborhoods $U\subset \mathcal L(V,E)$ and $W\subset \mathcal L(K,K^\perp)$ of $\alpha_0$ and $0$ respectively, and a smooth function $ \varphi : U\to W$  such that 
$$ \forall (\alpha, f)\in U\times W, \ F(\alpha, f)= 0\Leftrightarrow f=\varphi(\alpha).$$
Besides,
$ d\varphi(\alpha_0)=- \left(d_f F(\alpha_0,0)\right)^{(-1)} \circ d_\alpha F (\alpha_0,0),$
hence the result.
\epr
	Let $(V,g)$ and $(E,h)$ as in Lemma~\ref{ram}. Define
\beq\label{kappa}
\kappa: \mathcal L_{onto}(V,E)&\to& \grass(n-r,V)\\
\alpha& \mapsto & \ker \alpha.\nonumber
\eeq
The following lemma computes the derivative of $\kappa$.
	\begin{lemma}\label{ker}
	Let $(V,g)$ and $(E,h)$ be two real vector spaces as in Lemma~\ref{ram}.
	Then, $\kappa$ defined by~(\ref{kappa}) is smooth and for any $\alpha_0\in \mathcal L_{onto}(V,E),$
	in the chart given by Lemma~\ref{gra},
	\beqe 
	d\kappa(\alpha_0) : \mathcal L(V,E) &\to & 
	T_{\ker \alpha_0}\grass(n-r,V) \simeq 
	\mathcal L(\ker \alpha_0, \ker^\perp \alpha_0)\\
	\beta &\mapsto &-(\alpha_{0|K^\perp})^{-1}\beta_{|K}.
	\eeqe
\end{lemma}
\bpr
Let $\alpha_0\in \mathcal L_{onto}(V,E) $.
By Lemma~\ref{ram}, 
Locally, for any $\alpha$ close enough to $\alpha_0$, $$\kappa(\alpha) = (\id_{|\ker \alpha_0}+\varphi(\alpha))(\ker\alpha_0).$$
The second assertion of Lemma~\ref{ram} concludes. 
\epr
\subsection{The field and its geometry}

Let $n\geq 2$ and $1\leq r\leq n-1$ be integers, $(M,g)$ be a smooth Riemannian manifold of dimension $n\geq 1$, and $E\to M$ be a smooth 
real vector bundle of rank $r$. Define
$$ F= E\oplus T^*M\otimes  E.$$
Let $p : M\to \R$ be a Morse function, and $W\subset F$ be the subset of $F$ defined by
$$  \ 
W=\left\{\left(x,0,\alpha\right)\in F, \ dp(x)\neq 0,\  \alpha \text{ is onto and } \ker \alpha\subset \ker dp(x)\right\}.$$
Note that $W$ projects onto $M\setminus \crit(p).$ Since $p$ is Morse, $\crit(p)$ is a discrete set in $M$ without any accumulation point. If $M$ is compact, $\crit (p)$ is finite.
We will use later that for any $C^1$ section $s$ of $E$, then $x\in M$ is critical for the restriction of $p$ on $\{s=0\}$ at $x$ is equivalent  to $(x,s(x),\nabla s(x))\in W$, see~(\ref{x}) below. 
For any $x\in M\setminus \crit(p)$, let 
\beq\label{wx}
 W_x =\left\lbrace(0,\alpha)\in F_x, \ (x,0,\alpha)\in W\right\rbrace.
 \eeq

\begin{lemma}\label{sub} Let $M$, $E$, $p$, $F$ and $W$ be defined as above. Then,
	\begin{enumerate}
		\item \label{un}	$W$ is a  smooth submanifold of $F$ of codimension $n$;
		\item \label{deux} $W$ intersects transversally the fibres of $F$;
\item	\label{trois}	for any $x\in M\setminus \crit(p)$, $W_x$ is a smooth submanifold of $F_x$ of codimension $n$, and
		\beq\label{twx}
		 \forall (0,\alpha)\in W_x, \ T_{(0,\alpha)}W_x = \{(0,\beta)\in F_x, \  dp(x)(\alpha_{|\ker^\perp \alpha})^{-1}\beta_{|\ker \alpha} =0\}.
		 \eeq
	\end{enumerate}
	\end{lemma}
\bpr
Let $(x_0,0,\alpha_0)\in W$ and $K=\ker \alpha_0\subset T_{x_0}M$. Let $O$ be a neighborhood  of $x_0$   such that $M$ can be locally identified with $T_{x_0}M$ by a chart over $O$, and $E_{|O}$ can be identified with $O\times E_{x_0}$ via a trivalization. Then, $F_{|O}$ can be identified with $O\times (E_{x_0}\oplus T_{x_0}^* M \otimes E_{x_0}).$ 
By Lemma~\ref{ram}, there exists $U\subset \mathcal L(T_{x_0}M,E_{x_0})$ a neighborhood of $\alpha_0$ and a smooth map $\varphi : U \to \mathcal L (K,K^\perp)$, such that for any $\alpha\in U, $ $\ker \alpha  = (\id +\varphi(\alpha))_{|K}.$
Now, define the smooth map 
\beq\label{grosphi}
\Phi :  F_{|O} &\to &E_{x_0}\times K^*     \\
 (x,s,\alpha)&\mapsto &\left[s,dp(x)\left(\id +\varphi(\alpha)\right)_{|K}\right].\nonumber
\eeq
Then, $W\cap O = \Phi^{-1}(0).$
Moreover, by Lemma~\ref{ram} again, for all $ (v,t,\beta)\in T_{x_0 }M\times E_{x_0}\times \mathcal L(T_{x_0}M,E_{x_0}), $
\beq\label{diffphi}
 d\Phi (x_0,0,\alpha_0) (v,t,\beta)= \left[t,d^2p(x_0)(v)_{|K}-dp(x) (\alpha_{|K^\perp})^{-1}\beta_{|K}\right].
 \eeq
Since $dp(x)\neq 0$, $d\Phi$ is onto, so that $W$ is a smooth submanifold of $F$ of codimension $n$.
The third assertion of the lemma is an immediate consequence of~(\ref{diffphi}).
For the second assertion of the lemma, let $(x,0,\alpha)\in F\cap W $. Then, by~(\ref{diffphi}), 
$(v,t,\beta)\in F_x\cap T_{(x,0,\alpha)}W$ 
iff $(v,t)=0$ and 
$ dp(x) (\alpha_{|\ker^\perp \alpha})^{-1}\beta_{|\ker \alpha}=0,$ that is $(0,\beta)\in T_{(0,\alpha)}W_x$. Hence,
$F\pitchfork W$.
\epr

\subsection{Two Jacobians}

In this paragraph, the setting is the same as in the latter one, with the novelty that the vector bundle $E$ is endowed with a Euclidean metric $h_E$ on its fibres. We compute two Jacobians which will be needed for the coarea formula used in the main Kac-Rice formula Corollary~\ref{coro1}. 
For any  $(x,\alpha)\in TM^*\otimes E$, such that $\alpha$ is onto and 
$ dp(x)(\alpha_{|\ker^\perp\alpha})^{-1}\neq0\in E_x^*,$
 define
\beq\label{mumu}
\mu(x,\alpha)= \ker dp(x)(\alpha_{|\ker^\perp\alpha})^{-1}\subset E_x.
\eeq
Let $\epsilon(x,\alpha)$ be one of the two  unit vector in $\mu(x,\alpha)^\perp \subset E_x$, and $K=\ker\alpha$.
The following decomposition will help:
\beq\label{decompo}
 (\alpha_{|K^\perp})^{-1} = 
\begin{blockarray}{ccl}
  \R \epsilon & \mu & \\
  \begin{block}{(cc)l}
(\alpha_{|K^\perp})^{-1}_{|\R\epsilon} & 0 & \ker^\perp  dp(x)\\
*& (\alpha_{|K^\perp})^{-1}_{|\mu} &  \ker dp(x)\cap K^\perp \\
 \end{block}
 \end{blockarray}.
 \eeq
 Note that for any $x\in M\setminus \crit(p)$, by~(\ref{twx}), $$T_{(0,\alpha)}W_x=  
 \mathcal L(\ker \alpha, \mu(x,\alpha))\oplus
 \mathcal L(\ker^\perp \alpha,E_x)
 $$
 \begin{definition}\label{jx}(see~\cite[C.1]{stecconi}) Let $M,N$ be two Riemannian manifolds and $\kappa : M\to N$ be a $C^1$ map. Then, $J_x \kappa$ denotes the normal Jacobian, that is 
 the determinant in orthonormal basis of $d\kappa(x)_{|\ker^\perp d\kappa(x)}$. 
 \end{definition}
In the following, for any $x\in M$,  let $\kappa : \mathcal L_{onto}(T_xM,E_x) \to \grass(n-r,T_xM)$ defined by~(\ref{kappa}) for $V=T_xM$ and $E=E_x$. By an abuse of notation, we denote also $\kappa$ the  map $E_x\times \mathcal L_{onto}(T_xM,E_x)\ni (0,\alpha) \mapsto \kappa(\alpha).$
\begin{lemma}\label{ji} For any $x\in M\setminus \crit(p)$, 
	let $\kappa_{|W_x}: W_x\ni (0,\alpha)\mapsto \ker \alpha \subset \ker dp(x).$
	Then, for all  $(0,\alpha)\in W_x, \ J_{(0,\alpha)} (\kappa_{|W_x})	= \left|\det 
\alpha_{|\ker^\perp \alpha\cap \ker dp(x)}\right|^{n-r}$.
	\end{lemma}
\bpr
Firstly,  by Lemma~\ref{sub},
$$ T_{(0,\alpha)}W_x \ker(\beta\mapsto \langle \epsilon(x,\alpha),\beta_{|K}\rangle),$$
where $\epsilon$ has been defined above.
Since $$d(\kappa_{|W_x}) (0,\alpha)= (d\kappa(0,\alpha))_{|T_{(0,\alpha)} W_x},$$
from Lemma~\ref{ker} we infer that 
$ J_{(0,\alpha)} (\kappa_{|W_x})	= \left|\det 
	(\alpha_{|K^\perp})^{-1}_{|\mu(x,\alpha)}\right|^{n-r},$
	where $\mu$ has been defined by~(\ref{mumu}).
	Since $\alpha_{|K^\perp}$ induces an isomorphism between
	$ K^\perp \cap \ker dp(x)$ and $\mu(x,\alpha)$, we obtain the result. 
\epr
\begin{lemma}\label{g} Fix $x\in M\setminus \crit(p)$ and $K\in \grass(n-r,\ker dp(x))$. Let
	\beqe
	g:  \kappa^{-1}(K) & \to &  \grass(r-1,E_x)\\
	\alpha&\mapsto & \mu(x,\alpha)=
	\ker \left(dp(x)\circ (\alpha_{|K^\perp})^{-1}\right).
	\eeqe
	Then, for all $\alpha$, $ J_{(0,\alpha)} g = |\det 
	\alpha_{|\ker^\perp \alpha\cap \ker dp(x)}|$.
	\end{lemma}
\bpr 
Firstly, the map 
$$\nu: \alpha \in \kappa^{-1}(K)\mapsto (\alpha_{|K^\perp})^{-1}\in \mathcal L(E_x,K^\perp)$$ 
is smooth, and for any $\alpha \in \kappa^{-1}(K)$,
$$\forall \beta \in T_\alpha \kappa^{-1}(K), \  d\nu(\alpha)(\beta) = -
(\alpha_{|K^\perp})^{-1}\beta_{|K^\perp} (\alpha_{|K^\perp})^{-1}\in \mathcal L(E_x,K^\perp).$$
Moreover, by Lemma~\ref{kappa}, 
the differential of 
\beqe
\kappa_E : E_x^*\setminus\{0\}&\to &\grass(r-1,E_x)\\
f&\mapsto & \ker f
\eeqe
satisfies that, for any $f\in E_x^*\setminus\{0\}$ and $h\in E_x^*$, 
$ d\kappa_E(f)h= -(f_{|\ker^\perp f})^{-1}h_{|\ker f},$
	so that for any  $\alpha \in \kappa^{-1}(K)$ and any $\beta \in T_\alpha \kappa^{-1} K,$ 
	$$ dg(\alpha)\beta =
	\left(dp(x)(\alpha_{|K^\perp})^{-1}_{|\R\epsilon(x,\alpha)}\right)^{-1}
	dp(x)
	(\alpha_{|K^\perp})^{-1}\beta_{|K^\perp} (\alpha_{|K^\perp})^{-1}_{|\mu(x,\alpha)}\in \mathcal L(\mu(x,\alpha),\R \epsilon(x,\alpha)),
	$$
	hence
	$ J_{\alpha} g = |\det (\alpha_{|K^\perp})^{-1}_{|\mu(x,\alpha)}|.$
	Since $\alpha_{|K^\perp}$ induces an isomorphism between
$ K^\perp \cap \ker dp(x)$ and $\mu(x,\alpha)$, we obtain the result. 
\epr

\section{The mean number of induced critical points}\label{kr}

In the first part of this section, we provide two results. The first one, Proposition~\ref{kac}, is a Kac-Rice formula for the mean number of critical points of the restriction of a Morse function to the vanishing locus of a random section of some vector field. It is an application of the general Kac-Rice formula given by Theorem~\ref{stecconi}. The second result, Corollary~\ref{coro1}, is a more explicit and computable Kac-Rice formula which will be used in the applications of section~\ref{app}. In the second part of the section, we adapt the formula in a holomorphic context, see Corollary~\ref{thc}.
\subsection{The general formula}\label{formula}
Let $M$, $E$, $F$ and $W$ be as in section~\ref{prel}, and let $\nabla$ be a smooth connection on $E$. 
For any section $s\in C^1(M, E)$, let $ Z_s:= \{x\in M, s(x)=0\}$ and $X\in C^0(M,F)$ defined by
\beq\label{x}  X : M  & \to & F=E\oplus T^*M\otimes E \nonumber \\
x&\mapsto &\left[x,s(x), \nabla s(x)\right].
\eeq
Note that for any $x\in M$, 
$X(x)\in W$ if and only if $x\in Z_s$ and $x$ is a critical point of $p_{|Z_s}.$
For any random section $s\in  C^2(M,E) $, we are interested in the subset of $M$:
\beq \label{cris} \crit^p_i(s)=\{x\in M, \ X(x) \in W \text{ and } \Ind \nabla^2 (p_{|Z_s})(x)=i\},
\eeq
where $\nabla$ is any connection over $Z_s$. 
Note that $\Ind \nabla^2 (p_{|Z(s)})(x)$ is well defined for any $x\in X^{-1}(W)$ because in this case $d(p_{|Z_s})(x)=0$. Hence, $x\in X^{-1}(W)$ if and only if $p$ lies in $Z_s$ and the restriction of $f$ to $Z_s$ at $p$ is critical and its index equals $i$.
For any $s\in C^2(M,E)$ and $x\in X^{-1}(W)$, define also
\beq\label{pix}
\pi(x,\alpha):=\nabla^{2}p(x)_{|\ker \alpha} -
dp(x) \left(\alpha_{|\ker^\perp\alpha}\right)^{-1} \nabla^2 s(x)_{|\ker \alpha} \in \sym^2(\ker \nabla s(x),E),
\eeq
where $\alpha= \nabla s(x).$
Here, $\nabla^2 p$ denotes the covariant derivative of $dp$ for the Levi-Civita connexion associated to $g$. However, the formulas will not depend on the choice of this particular connexion. 
\begin{lemma}Assume that $p: M\to \R$ is a Morse function. Let  $s\in C^2(M,E)$ be a section of $E$, $x\in X^{-1}(W)$. Then,
	$$ \Ind(p_{|Z_s}(x)) = \Ind \left(\pi(x,\nabla s(x))\right).$$
\end{lemma}
\bpr Let $s\in C^2(M,E)$,  $x_0\in Z_{s}$ and  $K=\ker \nabla s(x_0)$. Assume that $\dim K=n-r$. 
We choose coordinates near $x_0$ so that $M$ is identified with $T_{x_0}M$.
By 
the implicit function theorem,
locally near $x_0=(0,0)\in K\oplus K^\perp$, and 
$Z_s$ is the graph over $K$ of a $C^2 $ map $f: K\to K^\perp$ with $f(0)=0$ and $df(0)=0$. 
Since locally $\forall z\in K, s(z,f(z))= 0$, we obtain
$$ \nabla_z s +\nabla_y s\circ df = 0, $$
where $\nabla_z$ and $\nabla_y$ denote the partial covariant derivatives along $K$ and $K^\perp$ respectively. 
so that$
\nabla^2_{z^2} s(0,0) +\nabla_y s(0,0) \circ d^2f(0) = 0, $
Now let $p_0: K\to \R$ define locally by 
$$ \forall z\in K, \ p_0(z)= p(z,f(z)).$$
Note that if $K\subset \ker dp(z)$, then 
$\Ind d^2 p(z) = \Ind \nabla^2 p_{|Z_s}$.
Now
$$dp_0= d_zp+d_y p\circ df,$$
so that
$ d^2 p_{0} (0) = 
d^2_zp(0,0)+d_y p\circ d^2f(0).$
Replacing $d^2f(0)$ by its value above, we obtain the result.
\epr

We will use the following general Kac-Rice formula. 
\begin{theorem}(\cite[Theorem 3.3]{stecconi}\label{stecconi}) Let $n$ be a positive integer, 
	$M$ be a smooth manifold of dimension $n$, 
$F \to M$ be a smooth vector bundle and $X \in \Gamma(M,F)$ be a non-degenerate smooth Gaussian random section. Let $W \subset F$ be a smooth submanifold of codimension $n$ such that for every $x \in  M$, $W_x:=W\pitchfork F_x$. Let the total space of $F$ be endowed
	with a Riemannian metric that is Euclidean on fibers. Then for any Borel subset $A \subset M$
	\beqe
		\mathbb E \# \{x\in A\cap X^{-1}(W)\} = \int_{x\in A} 
	\int_{q\in W_x}& \mathbb E \left( 
	J_{x} X\frac{\sigma_q (X,W )}{\sigma_q(F_x,W)} | X(x) = q\right) 
	\rho_{X(x)}(q)d\vol(q) d\vol(x),
	\eeqe
	where 
$	 \rho_{X(x)}(q)	$ is the density of $X(x)$ at $q$
	and besides, $\sigma_q (X,W )$, $\sigma_q(F_x,W)$ denote the “angles”
	made by $T_qW$ with, respectively, $d_x X(T_xM)$ and $T_qF_x$, see~\cite[Definition B.2]{stecconi}.
\end{theorem}
The random section $X \in \Gamma(M,F)$ is said to be \emph{non-degenerate}~\cite[Definition 3.1]{stecconi} if  for 
any $x\in M$, $\text{supp }X(x) = F_x$.
We will not explain here the terms $\sigma_q$, because by the proof of \cite[Lemma 7.2]{stecconi}, locally
\beq\label{delta}
 J_{x} X\frac{\sigma_q (X,W )}{\sigma_q(F_x,W)} =
\frac{\delta_x (\Phi \circ X)}{J_x (\Phi_{|F_x})},
\eeq
where $\Phi : F\to \R^n$ is a local defining 
function for $W$, that is $W=\Phi^{-1}(0)$, where 
$ \delta_p $ denotes the Jacobian density~\cite[(A.1)]{stecconi} and where 
$J_x$ is the normal Jacobian, see Definition~\ref{jx}.

The following Proposition~\ref{kac} is an application of the general Kac-Rice formula above, namely a Kac-Rice formula 
for the number of induced critical points of the restriction of a Morse function on random nodal sets. We need some notations. Let $s\in C^2(M,F)$, $x\in X^{-1}(W)$ and $\alpha = \nabla s(x)$.  Recall that $\epsilon(x,\alpha)\in E_x$ denotes a unit vector of $\ker^\perp dp(x) (\alpha_{|K^\perp})^{-1}\subset E_x$, that  $\pi(x,\alpha)$ is defined by~(\ref{pix}) and $\crit^p_i$ by~(\ref{cris}). Also, let $h(x)$ be one of the two unit vector in $\ker^\perp dp(x)\subset T_x M$. 
\begin{proposition}\label{kac}Let $n\geq 2$ and $1\leq r\leq n-1$ be integers, $(M,g)$ be a Riemannian manifold, $(E,h)\to M$ be a rank $r$ smooth Euclidean vector bundle and $s\in \Gamma(M,E)$ be a non-degenerate Gaussian smooth field. Let $p: M\to \R$ be a smooth Morse function.
Then, for any $i\in \{0,\cdots,n-r\}$
and any	Borel subset $A\subset M$, 
\beqe	\mathbb E \left[\# \left(A\cap \crit_i^p\right)\right] &=& \int_{x\in A}
\int_{\substack{\alpha\in \mathcal L_{onto}(T_xM,E_x)\\ \ker \alpha\subset \ker dp(x)}} 
\left|\det \alpha_{|\ker^\perp \alpha}	\right| \\
&&\mathbb E \Big[	{\bf 1}_{\{\Ind (\pi(x,\alpha))=i\}} 
\Big| \det \Big(
\langle \nabla^2 s(x)_{|\ker \alpha},\epsilon(x,\alpha)\rangle
\\ &&
- \langle \alpha(h(x)),\epsilon(x,\alpha)\rangle \frac{\nabla^{2}p(x)_{|\ker \alpha}}
{\|dp(x)\|}
\Big)\Big|
 \left. \middle| \  s(x)=0, \nabla s(x)=\alpha\right.\Big]
 \\ &&
\rho_{X(x)}(0,\alpha)
d\vol(\alpha)  d\vol(x),
\eeqe
	where $\rho_{X(x)}$ is the Gaussian density of $X(x)$
	and the determinants are computed in orthonormal basis.
	Moreover, this integral is finite if $\vol (A)$ is finite.
\end{proposition}
\bpr We use Theorem~\ref{stecconi}, using locally~(\ref{delta}). So let $(x_0,0,\alpha_0) \in W$.
Locally and in coordinates, 
using the local defining function $\Phi$ for $W$ given by~(\ref{grosphi}),
$$ \forall x\in O, \ \Phi (X(x))= 
[s(x), dp(x)(\id +\varphi(\nabla s(x))]\in E_{x_0}\times (\ker \alpha_0)^*,$$
where $\nabla$ still denotes the connection $\nabla$ through the trivialization. 
Hence,
\beqe
\forall v\in T_xM,\  d(\Phi\circ X)(x)(v)&=& 
\Big[ds(x)(v), d^2p(x)(v)(\id +\varphi(\nabla s(x)))_{|\ker\alpha_0}\\
&& -dp(x)(\nabla s(x)_{|\ker \nabla 	f(x)^\perp})^{-1} d\nabla s(x)(v)_{|\ker \alpha_0} \Big].\nonumber
\eeqe
In particular, for any $x\in X^{-1}(W)$ and any $v\in T_x M$, if $K=\ker \nabla s(x)$, 
$$ d(\Phi\circ X)(x)(v)= 
\Big[\nabla_v s(x), \nabla_v dp(x)_{|K}-dp(x)(\nabla s(x)_{|K^\perp })^{-1} \nabla_v \nabla s(x)_{|K}\Big].$$
Now, decomposing $T_xM$ as
$T_xM= K^\perp\oplus K$, since $\nabla s(x)_{|K}=0
$,
computing the determinant of this differential gives
$$ \delta_x(\Phi\circ X)= 
\left|\det \nabla s(x)_{|K^\perp}\right| \left|\det \left(\nabla  dp(x)_{|K}-dp(x)(\nabla s(x)_{|K^\perp })^{-1}  \nabla^2 s(x)_{|K}\right)\right|.$$
Now, let us compute $J_x \Phi_{|F_x}$. For this, recall that
$$ \forall (s,\alpha)\in F_x,\ 
\Phi(s,\alpha)=
[s, dp(x)(\id +\varphi(\alpha))],$$
so that 
$$ \forall (0,\alpha)\in W_x,\ \forall (t,\beta)\in T_pF_x,\ 
d(\Phi_{|F_x})(0,\alpha)(t,\beta)=
[t, -dp(x)(\alpha_{|K^\perp})^{-1}\beta_{|K}].$$
Since 
$ J_{(0,\alpha)}
(\beta \mapsto \langle \beta_{|K}, \epsilon(x,\alpha)\rangle)= 1,$
we get that 
$$ J_{(0,\alpha)} (\Phi_{|F_x })
= |dp(x)(\alpha_{|K^\perp})^{-1}\epsilon(x,\alpha)|^{n-r}.$$
Moreover,
\beqe
dp(x) (\alpha_{|K^\perp})^{-1} \nabla^2 s(x)_{|K}
= \left(dp(x)(\alpha_{|K^\perp})^{-1}\epsilon(x,\alpha)\right) \langle\nabla^2 s(x)_{|K},\epsilon(x,\alpha)\rangle.
\eeqe
By Lemma~\ref{sub}, $W$ intersects the fibres of $F$ transversally, so that Theorem~\ref{stecconi} applies. Replacing the integrant in the theorem by~(\ref{delta}), we obtain the formula.
Finally, by the decomposition~(\ref{decompo}),
$$ dp(x) (\alpha_{|\ker^\perp \alpha})^{-1}\epsilon(x,\alpha) = \|dp(x)\| (\langle \alpha(h(x)),\epsilon(x,\alpha)\rangle)^{-1}.$$
Since $p$ is Morse, for any critical point $x\in \crit(p)$, there exists a constant $C_x$ such that  $\|dp(y)\|\geq C_x \|y-x\|.$ Hence, the pole in the integration over $U$ created by $x$ has order $n-r$, which is integrable, see also~\cite[Remark 3.3.3]{gayet2015expected}.
\epr

In order to provide an effective formula in concrete settings, we add further parameters.
For any $x\in M$,  recall that $h(x)\in T_xM$ be one of the two unit vectors in $\ker^\perp dp(x)$. Moreover, for any real hyperplane $\mu\subset E_x$, let $\epsilon(\mu)$ be a unit vector in $\mu^\perp$. Recall that $\crit_i^p$ is defined by~(\ref{cris}).
\begin{corollary}\label{coro1}
Assume the hypotheses of Proposition~\ref{kac} are satisfied.  
Let	$A\subset M$ be a Borel subset. Then, for any $i\in \{0,\cdots,n-r\},$
	\beqe	\mathbb E \left[\# \left(A\cap \crit_i^p\right)\right] = 
	\int_{x\in A}&&
	\int_{\substack{K\in \grass(n-r, \ker dp(x))\\
			\mu\in \grass(r-1,E_x)}} 
	\int_{\substack{\alpha\in TM_x^*\otimes E_x
			\\ \ker \alpha=K\\
			\ker dp(x)(\alpha_{|K^\perp})^{-1}=\mu	
	}} \\
	&&	\left|\det (\alpha_{|K^\perp\cap \ker dp(x)})\right|^{n-r+2} 
	|\langle \alpha(h(x)),\epsilon(\mu)\rangle|\\
	&&\mathbb E \Big[	{\bf 1}_{\{\Ind (\pi(x,\alpha))=i\}} 
	\Big| \det \Big(
		\langle \nabla^2 s(x)_{|K},\epsilon(\mu)\rangle
		\\ &&-
		 \langle \alpha(h(x)),\epsilon(\mu)\rangle \frac{\nabla^{2}p(x)_{|K}}
		{\|dp(x)\|}
	\Big)\Big|
	 \left. \middle| \  s(x)=0, \nabla s(x)=\alpha\right.\Big]
	\\ && \rho_{X(x)}(0,\alpha)
	d\vol(\alpha)d\vol(\mu)	d\vol(K)  d\vol(x),
	\eeqe
	where $\pi(x,\alpha)$ is given by~(\ref{pix}) and $\rho_X$ denotes the density of $X$. 
	\end{corollary}
\bpr In the formula given by Proposition~\ref{kac},
we handle first the determinant of $\alpha_{|K^\perp}$. 
Since $(\alpha_{|K^\perp})^{-1}(\mu)=\ker dp(x)\cap K^\perp,$
if $h(x)$ is a unit vector in $\ker^\perp dp(x)$, then
\beq\label{bernin}
 | \det \alpha _{|K^\perp}| = 
|\det \alpha_{|K^\perp\cap \ker dp(x)}| |\langle\alpha (h(x)),\epsilon(\mu)\rangle|.
\eeq
We then apply two times the coarea formula (see for instance~\cite[Theorem C.3]{stecconi} from which we borrow the notations) for the integral in $\alpha$. The first formula is applied with the map $\kappa_{|W_x} : W_x \to \grass(n-r, \ker dp(x))$, where $\kappa$ is  defined by~(\ref{ji}). By Lemma~\ref{ji}, its Jacobian satisfies, for any $(0,\alpha)\in W_x$, $J_{(0,\alpha)} (\kappa_{|W_x})	= \left|\det 
\alpha_{|\ker^\perp \alpha\cap \ker dp(x)}\right|^{n-r}$.
The second coarea formula is applied with $K\in \grass(n-r,\ker dp(x))$ fixed, with the function 
$g:  \kappa^{-1}(K) \to  \grass(r-1,E_x)$ defined in Lemma~\ref{g}. Then, By the latter, for all $\alpha$, $ J_{(0,\alpha)} g = |\det 
\alpha_{|\ker^\perp \alpha\cap \ker dp(x)}|$. We obtain the result. Together with~(\ref{bernin}), we obtain the desired formula.
\epr

\subsection{The holomorphic setting}

In this paragraph, let $n\geq 2$ and $1\leq r\leq n-1$ be  integers, $M$ be a complex smooth manifold of complex dimension $n$, endowed with a Hermitian metric $g$. Let $(E,h_E)\to M$ be a holomorphic Hermitian vector bundle of rank $r$, and $s \in \Gamma(M,E)$ be a holomorphic Gaussian field. In section~\ref{app}, $M$ will be either a compact projective manifold and $E$ the tensor product of a fixed vector bundle tensored by the high powers of an ample line bundle, or $M$ will be the affine complex space and $E$ the trivial complex vector bundle of rank $r$. Let $\nabla $ be the Chern connection for $E$, that is the unique holomorphic and metric connection on $E$, see~\cite{griffiths}. In this complex case, the real setting of paragraph~\ref{formula} adapts formally, changing the field $\R$ into $\C$. 
In particular, we define
$$ F= E\oplus \mathcal L^\C (TM,E).$$
However specific changes must be also done. Let $p : M\to \R$ be a smooth Morse function. Then, for any holomorphic section $s$ of $E$ and any $x\in Z_s$, 
$$\ker \nabla s(x)\subset \ker dp(x)\Leftrightarrow
\ker \nabla s(x) \subset \ker \pi_\C(x),$$
where $\pi_\C(x)$ denotes the complexification
of $dp(x)$, that is $\pi_\C(x)\in \mathcal L^\C (T_xM, \C)$
and $dp(x) = \Re \pi_\C(x) $. Then, 
we use that for any complex subspace $K\subset T_x M$ and any $\alpha\in \mathcal L^\C(K,E_x)$, 
the real determinant (computed in orthonormal basis) of the associated real map $\alpha_\R$ equals
\beq \label{dety}
|\det \alpha_\R| = |\det \alpha|^2.
\eeq
As in the real case, the Gaussian holomorphic field $s$ is said to be \emph{non-degenerate} if for any $x\in M$, 
$$s\mapsto (s(x),\nabla s(x))\in E\times \mathcal L^\C (T_xM, E_x)$$
is onto. 
As before, we define $$W=\{(x,0,\alpha)\in F, dp(x)\neq 0, \alpha \text{ onto and } \ker \alpha\subset \ker \pi_\C (x) \},$$
and $W_x$ as the fibre of $W$ over $x$.
For any $K\in \grass_\C(n-r, T_x M)$ and $\mu\in \grass_\C(r-1, E_x)$, let
$$ W(x,K,\mu):=\{\alpha\in \mathcal L_{onto}^\C(TM_x,\otimes E_x) \, |\, 
\ker \alpha=K, \ 
\ker \pi_\C(\alpha_{|K^\perp})^{-1}=\mu	
\}.$$
\begin{lemma}\label{doublev} Under the hypotheses above, $W(x,K,\mu)$ 	is a submanifold of $W_x$ of complex dimension
$r^2 - (r-1).$
	\end{lemma}
For any $s\in H^0(M,E)$ and $x\in X^{-1}(W)$, define also
\beq\label{pixc}
\pi(x,\alpha):=\nabla^{2}\pi_\C(x)_{|\ker \alpha} -
\pi_\C(x) \left(\alpha_{|\ker^\perp\alpha}\right)^{-1} \nabla^2 s(x)_{|\ker \alpha} \in \sym^2(\ker \nabla s(x),E_x),
\eeq
where $\alpha= \nabla s(x).$
Lastly,
for any $x\in M$, denote by $h(x)\in T_xM$  any  unit vector in $(\ker \pi_\C(x))^\perp\subset T_x M$,
and for any complex hyperplane $\mu \subset E_x$, let $\epsilon(\mu)$ be a unit vector in $\mu^\perp\subset E_x$.
Recall that $\crit_i^p$ is defined by~(\ref{cris}).
\begin{theorem}\label{thc}
Let	$(M,g)$ be a complex manifold, $(E,h_E)\to M$ be a holomorphic Hermitian vector bundle, and $s \in \Gamma(M,E)$ be a non-degenerate holomorphic Gaussian field.
Let $A\subset M$ any Borel subset. 
Then, 
\beqe	\mathbb E \left[\# \left(A\cap \crit_i^p\right)\right] = 
\int_{x\in A}&&
\int_{\substack{K \in \grass_\C(n-r, \ker \pi_\C(x))\\
		\mu\in \grass_\C(r-1,E_x)}} 
\int_{\substack{\alpha\in \mathcal L^\C(T_xM, E_x)
		\\ \ker \alpha=K\\
		\ker \pi_\C(x)(\alpha_{|K^\perp})^{-1}=\mu	
}} \\
&&	\left|\det (\alpha_{|K^\perp\cap \ker \pi_\C(x)})\right|^{2(n-r+2)} |\langle\alpha (h),\epsilon(\mu)\rangle|^2 \\
&&\mathbb E \Big[	{\bf 1}_{\{\Ind (\pi(x,\alpha))=i\}} 
\Big| \det{}_\R \Big(
\langle \nabla^2 s(x)_{|K},\epsilon(\mu)\rangle\\
&& -
		 \langle \alpha(h(x)),\epsilon(\mu)\rangle \frac{\nabla \pi_\C(x)_{|K}}
{\|\pi_\C(x)\|}
\Big)\Big|
 \left. \middle| \  s(x)=0, \nabla s(x)=\alpha\right.\Big]\\
 &&
\rho_{X(x)}(0,\alpha)
d\vol(\alpha)d\vol(\mu)	d\vol(K)  d\vol(x),
\eeqe
where $\pi(x,\alpha)$ is given by~(\ref{pixc}) and $\rho_X$ is the density of $X$. Moreover, the integral is finite if $\vol (A)$ is finite. 
\end{theorem}
\bpr The proof is formally the same as the one of Corollary~\ref{coro1}, using the rules mentionned above, 
so we omit it.
\epr

\section{Applications}\label{app}
In this section we apply Theorem~\ref{thc} to the complex Bargmann-Fock field on $\C^n$ and then to the projective setting. Finally, we apply Proposition~\ref{kac} to the boundary case, which is a mixed between complex and the real setting and is needed for the main Theorems~\ref{main} and~\ref{BF}.

\subsection{The Bargmann-Fock field }

Recall that the	Bargmann-Fock field is defined
by 
\beq\label{bafo}
 \forall z\in \C^n, \ f(z)=\sum_{(i_1,\cdots, i_n)\in \Nn^n}a_{i_0, \cdots, i_n}
{\sqrt{\frac{\pi^{i_1+\cdots +i_n}}{i_1!\cdots i_n!}}}{z_1^{i_1}\cdots z_n^{i_n}}e^{-\frac12 \pi\|z\|^2},
\eeq
where the $a_I$'s are independent normal complex Gaussian random variables. The associated covariant function equals
\beq\label{mp}
\forall z,w\in \C^n, \ \mathcal P(z,w):= \mathbb E (f(z)\overline{f(w)})= \exp\left(
-\frac{\pi}2 (\|z\|^2 +\|w\|^2 -2\langle z,w\rangle_{\C^n})
\right).
\eeq
Even if the kernel $\mathcal P$ is not invariant under translation or rotations, the law of $Z_f$ is, see~\cite[Proposition 2.3.4]{hough2009zeros}.
The a priori superfluous presence of $\pi$ is in fact consistent with the projective situation. Indeed, the affine Bargmann-Fock is 
the universal local limit of the projective model, see Theorem~\ref{Dai}.
In order to unify the setting, 
we consider here that $M=\C^n$ and $L=\C^n \times \C$ with its standard Hermitian metric. 
Then $$\mathcal P (z,w)\in L_z\otimes L_w^*.$$
Note that the trivial connection on $L$ has vanishing curvature. Hence, let $\nabla_0$ be the metric connection for this setting:
\beq\label{nabla}
\nabla_0 1 = \frac12 \pi (\bar \partial -\partial) \|z\|^2,
\eeq
whereas the dual connection $\nabla^*$ on $L^*$
satisfies 
\beq\label{nabla2}
\nabla_0^* 1^* = -\frac12 \pi (\bar \partial -\partial) \|z\|^2,
\eeq
where $1^*$ is the dual of $1$. 
Notice that $1$ is no longer a holomorphic section for this connection, but 
the (peak) section (see~\cite{tian1990set}, \cite{donaldson1996symplectic}) $$\sigma_0:= \exp(-\frac12 \pi \|z\|^2)$$ is, 
since 
$ \nabla_0^{(0,1)}\sigma_0= (-\frac{1}2\pi \bar \partial \|z\|^2 +\frac{1}2\pi\bar \partial \|z\|^2) \sigma_0=0.$
The connection $\nabla_0$ is then the Chern connection for the trivial metric and this holomorphic structure. 
This implies that the section $\mathcal P$ is holomorphic in $z$, and antiholomorphic in $w$. 
Moreover, the
curvature of $\nabla_0$ equals
$$ \mathcal R_0 = \bar\partial \partial \log \|\sigma_0\|^2 = \pi \partial \bar \partial \|z\|^2,$$
and the
curvature form equals
$$\frac{i}{2\pi}\mathcal R_0=\frac{i}2 \sum_{i=1}^n dz_i \wedge \overline{dz_i}
$$
which is the standard symplectic form $\omega_0$ over $\R^{2n}$.
Now, almost surely an instance $f$ of the  Bargmann Fock Gaussian field is a holomorphic section
for the standard complex structure and the connection with standard curvature form given by~(\ref{nabla}).

Let $E=\C^n \times \C^r$ endowed with its trivial metric and 
let $f=(f_i)_{i=1, \cdots, r}$ be $r$ 
independent copies of the Bargmann-Fock field.  Then,
$f$ is a random section of $E\otimes L$,
and its covariance function equals $\mathcal P \text{Id}_{\C^r}.$
In the following theorem, recall that $\crit^p_i$ is defined by~(\ref{cris}), where we use the connexion $\nabla_0^r$ acting on sections of $E\otimes L$. By an abuse of notation, we continue to use $\nabla_0$ for $\nabla_0^r$. 
\begin{theorem}\label{BF0}Let $1\leq r\leq n $ be integers, $ f : \C^n\to \C^r$ be $r$ independent copies of the Bargmann-Fock field~(\ref{bafo}), and $U\subset \C^n$ be an open subset of finite volume. Let $p : U\to \R$ be a smooth Morse function. Then,
	\beqe \forall 0\leq i\leq 2n-2r\setminus\{n-r\},\ 
	\ \frac{1}{R^{2n}}\mathbb E \#(RU\cap \crit_i^p )&\underset{R\to +\infty}\to& 0 \\
	\frac{1}{R^{2n}}	\mathbb E \#(RU\cap \crit^p_{n-r} )&\underset{R\to +\infty}\to & n!{n-1 \choose r-1}\vol(U),
	\eeqe
	where $\vol $ denotes the volume for the standard metric on $\C^n$.
\end{theorem}
We postpone the proof of this theorem after the projective case, since the latter is similar but more complicated. In both cases, we will need 
the following lemma:
\begin{lemma}\label{deriv}Let $\mathcal P$ the Bargmann-Fock covariance~(\ref{mp}), and $\nabla_0$ the connection defined by~(\ref{nabla}) and~(\ref{nabla2}). Then, for any $z\in \C^n$, 
	\beqe \nabla_{0\  z,\bar w}^{(1,0),(0,1)}\mathcal P(z,z) &=& 
	\pi \sum_{i=1}^n dz_i\otimes \overline{dw_i}\\
\text{	and }
	 \nabla_{0 \ z^2,\bar w^2}^{(1,0)^2,(0,1)^2}\mathcal P(z,z) &=& 
	\pi^2 \sum_{i,j,k,\ell=1}^n(\delta_{ik}\delta{j\ell}+\delta_{i\ell}\delta_{jk}) dz_i \otimes dz_j\otimes \overline{dw_k}\otimes \overline{dw_\ell}.
\eeqe
\end{lemma}
\bpr This is a straightforward consequence of the definition of $\nabla_0$ and $\mathcal P$.\epr
We will need the following covariance matrix for Hessians: 
\beq\label{goe}
 \Sigma_{GOE} =  \left(\delta_{(ij)(kl)}+\delta_{(ji)(kl)}\right)_{\substack{1\leq i\leq j\leq n\\
		1\leq k\leq l\leq n}}\in M_{\frac{n(n+1)}2}(\C).
\eeq
\begin{corollary}\label{coaff} Let $f:\C^n \to \C^r$ be $r$ independent copies of the Bargmann-Fock field. 
	Then, for any $x\in \C^n,$
	\beqe
	\cov (f(x), \nabla_0 f(x),\nabla_0^2f(x)) &=& \begin{pmatrix}
		1 & 0 & 0\\
		0   & \pi \ \id_{\C^n}  & 0\\
		0&0& \pi^2\Sigma_{GOE}
	\end{pmatrix} \otimes \ \id_{\C^r}.
	\eeqe
\end{corollary}
\bpr This is an immediate consequence of Lemma~\ref{deriv}.\epr

 \subsection{The complex projective case}\label{projo}

Let $n\geq 2$, $1\leq r\leq n-1$ be integers, $M$ be a compact smooth complex manifold of dimension $n$ equipped with a holomorphic Hermitian vector bundle $(E,h_E)$ of rank $r$ and an ample holomorphic line bundle $(L,h)$. Assume that $h$ has a positive curvature form $\omega$, see~(\ref{omega}). Let $\nabla$ be the Chern connection of $E\otimes L^d$. 
Recall that $\crit_i^p$ is defined by~(\ref{cris}).
\begin{theorem} \label{projective}
	Let $M$, $(E,h_E)$, $(L,h)$, $\omega$ as above and let	$U\subset M$ be a 
	$0$-codimension submanifold with finite volume. Then
	\begin{eqnarray*}
		\forall i\in \{0,\cdots, 2n-2r\}\setminus\{n-r\}, \ \frac{1}{d^n}\mathbb E \#(U\cap \crit_i^p) &\underset{d\to \infty}{\to}& 0\\
		\frac{1}{d^n}\mathbb E \#(U\cap \crit_{n-r}^p)&\underset{d\to \infty}{\to} & {n-1 \choose r-1}\int_U \omega^n.
	\end{eqnarray*}	
The probability measure $\mu_d$ used for the average is defined by~(\ref{mesure}).
\end{theorem} 
Theorem~\ref{gw} (\cite[Theorem 3.5.1]{gayet2015expected}) implies 
this result for the squared modulus of a Lefschetz pencil $p: M\dashrightarrow  \C P^1$. Indeed,
since $p$ is holomorphic (outside its singular locus), 
$p_{|Z_s}$ is critical if and only if $|p|^2_{|Z_s}$ is, 
and in the latter case the index equals $n-r$.

\paragraph{Bergman and Bargmann-Fock. }
The covariance function for the Gaussian field generated by the holomorphic sections $s\in H^0(M,E\otimes L^d)$ is
$$ \forall z,w\in M, 
\ E_d (z,w) = \mathbb E \left[s(z)\otimes (s(w))^*\right]\in (E\otimes L^d)_z\otimes (E\otimes L^d)_w^*,$$
where $E^*$ is the (complex) dual of $E$ and 
$$\forall w\in M, \ \forall s,t\in (E\otimes L^d)_w, \  s^* (t) = h_E\otimes h_{L^d}(s,t).$$
The covariance $E_d$ is the \emph{Bergman kernel}, that is the kernel of the orthogonal projector from $L^2(M,E\otimes L^d)$ onto $H^0(M,E\otimes L^d)$. This fact can be seen through the equations
$$\forall z,w\in M, \ E_d(z,w) = \sum_{i=1}^{N_d} S_i (z)\otimes S_i^*(w),$$
where $(S_i)_i$ is an orthonormal basis of $H^0(M,E\otimes L^d)$ for the Hermitian product~(\ref{herm}).
Recall that the metric $g$ is induced by the curvature form $\omega$  and the complex structure. 
It is now classical  that the Bergman kernel has a universal rescaled (at scale $1/\sqrt d$) limit, the Bargmann-Fock kernel $\mathcal P$ defined by~(\ref{mp}). Theorem~\ref{Dai} below quantifies this phenomenon. For this, we need to introduce local trivializations and charts. 
Let $x\in M$ and $R>0$ such that $2R$ is less than the radius of injectivity of $M$ at $x$. Then the exponential map based at $x$ induces a chart near $x$ with values in $B_{T_xM}(0,2R)$. The parallel transport provides a trivialization 
$$\varphi_x : B_{T_xM}(0,2R)\times (E\otimes L_d)_{x}\to (E\otimes L_d)_{|B_{T_xM}(0,2R) }$$ 
which induces a trivialization of 
$(E\otimes L_d)\boxtimes
(E\otimes L_d)^*_{|B_{T_xM}(0,2R)^2}$.
Under this trivialization, the Bergman kernel 
$ E_d$ becomes a map from $T_x M^2 $ with values into $\End\left((E\otimes L^d)_x\right)$.

\begin{theorem}(\cite[Theorem 1]{ma2013remark})\label{Dai} Under the hypotheses of Theorem~\ref{main}, let $m\in \Nn$. Then, there exist $C>0$, such that for any $k\in\{0,\cdots, m\}, $ for any $x\in M$, 
	$\forall z,w\in B_{T_xM}(0,\frac{1}{\sqrt d}), $ 
	\beqe \left\|	D^k_{(z,w)}\left(
	\frac{1}{d^n}E_d(z,w) - \mathcal P(z\sqrt d,w\sqrt d )
	\ \id_{(E\otimes L^d)_x}
	\right)
	\right\|	
	\leq C d^{\frac{k}2-1}.
	\eeqe
\end{theorem}
The original reference is a little more intricated, see~\cite[Proposition 3.4]{letpuch} for the present simplification.
We will also need the following lemma:
\begin{lemma}\label{connecte}
	Under the local trivializations given before, at $x$ (the center of the chart) the two equalities hold: 
	$$
	\nabla = \nabla_0 +O(\frac{1}{\sqrt d})  \text{ and } \nabla^2 = (\nabla_0)^2
	+O(\frac{1}{\sqrt d}) .
$$
	\end{lemma}
\bpr The conjonction of \cite[Lemma 1.6.6]{ma2007holomorphic} and \cite[(4.1.103)]{ma2007holomorphic} implies 
that
$$\nabla = \nabla_0 +O(\frac{1}{\sqrt d}) +O(\|z-x\|^3), $$ which gives the first estimate. The second one is implied by the first one and by the fact that the Levi-Civita connection associated to $g$ is trivial at $x$, because the coordinates on $M$ are normal at $x$.
\epr
\begin{corollary}\label{coco} Under the hypotheses and trivializations above near $x\in M$, 
	in any orthonormal basis of $T_xM$, 
	\beqe \cov\left(s, \nabla s, \nabla^2 s\right)_{|x} =
	d^n\left(	\begin{array}{ccc}
		(1+O(\frac1d))& O(\frac{1}{\sqrt d}) & O(1)\\
		O(\frac{1}{\sqrt d}) & \pi d{I}_{n}(1+O(\frac1d))&  O(\sqrt d)\\
		O(1)& O(\sqrt d) &  \pi^2 d^2\Sigma_{GOE}(1+O(\frac1d))
	\end{array}\right) \id_{(E\otimes L^d)_x},
	\eeqe
	where $I_n\in M_n(\R)$ and $\Sigma_{GOE}$ is defined 
	by~(\ref{goe}).
	Moreover, for any $\alpha\in T^*_xM\otimes E_x$, 
 $$\left(\langle \nabla^2 s,\epsilon\rangle \ | \  s=0, \nabla s=\alpha\right) \sim N\left(O\big(\frac{\|\alpha\|}{\sqrt d}\big), \Sigma\right),$$
where
$$\Sigma :=\pi^2d^{n+2}\Sigma_{GOE}  \ \id_{(E\otimes L^d)_x}\big(1+O(\frac1d)\big).$$ The constants involved in the error terms do not depend on $\alpha$. 
	\end{corollary}
\bpr The first assertion is a direct consequence of Theorem~\ref{Dai}, Lemma~\ref{connecte} and Corollary~\ref{coaff}. The second one is deduced from the classical regression formula and from 
\beq (\cov(s, \nabla s))^{-1} =
\frac{1}{d^{n}}\left(	\begin{array}{cc}
	(1+0(\frac1d))& O(\frac{1}{d^{\frac32}}) \\
	O(\frac{1}{d^{\frac32}}) & \frac{1}{\pi d}(1+O(\frac1d))
\end{array}\right) \id_{(E\otimes L^d)_z}.
\eeq
\epr

\bpr[ of Theorem~\ref{projective}] We want to apply Theorem~\ref{thc}.
First, from Corollary~\ref{coco} we get that 
for any $x\in M$ and any $\alpha \in \mathcal L^\C(T_xM, E_x),$
\beqe 
\rho_{X(x)}(0,\alpha) = \frac{(1+O(\frac1d))}{(2\pi)^{r+nr}(d^n)^r(\pi d^{n+1})^{nr}}
\exp\left(-\frac{1}2 \frac{1}{\pi d^{n+1}}(1+O(\frac1d))\|\alpha\|^2
\right).
\eeqe
Now, if $K=\ker \alpha$, $\mu= \ker dp(x)(\alpha_{|K^\perp})^{-1}$ and $\epsilon (\mu)\in \mu^\perp $ has a norm equal to 1, let 
	Let  
	\beq\label{bab}
	(\beta,a,b) = \frac{1}{\sqrt{\pi d^{n+1}}}\left(
	\alpha_{|K^\perp_\C\cap \ker \pi_\C},
	\langle \alpha_{|\ker^\perp \pi_\C}, \epsilon\rangle,
	\pi^\perp_{\mu}\alpha_{|\ker^\perp \pi_\C}\right),
\eeq
	where $\pi^\perp_{\mu}$ denotes the orthgonal projection (for $h_E$) onto $\mu$. 
Using Lemma~\ref{doublev} for the transformation of $d\vol(\alpha)$, the term
	$$ 	\left|\det \alpha_{|K^\perp_\C\cap \ker \pi_\C}\right|^{2(n-r+2)} |\langle\alpha (h(x)),\epsilon\rangle|^2 d\vol (\alpha)_{|W(x,K,\mu) }
	\rho_{X(x)}(0,\alpha)
		$$ 
		in the integral of Theorem~\ref{thc} 
equals
\beqe
(1+O(\frac1d)) \frac{(2\pi)^{r-1}(\pi d^{n+1})^{(r-1)(n-r+2)+1+r^2-(r-1)}(2\pi)^{(r-1)^2+1}}
{(2\pi)^{r+nr}d^{nr}(\pi d^{n+1})^{nr}}\\
	\frac{|\det \beta|^{2(n-r+2)}}{(2\pi)^{(r-1)^2}} \frac{|a|^2}{2\pi} \frac{d\vol (\beta,a,b)}{(2\pi)^{r-1}}
\exp\left(-\frac12(1+O(\frac1d)) (\|\beta\|^2+|a|^2+\|b\|^2)\right).
\eeqe
Note that $$ \int_{a\in \C }|a|^2 e^{-\frac12|a|^2|}da = 4\pi.$$
By Corollary~\ref{coco}, the field $X$ defined by~(\ref{x}) is non-degenerate for $d$ large enough. Hence, we can apply Theorem~\ref{thc}.
Let $Y=\frac{1}{\sqrt	{\pi^2d^{n+2}}} \nabla^2 s(x).$
Then,
the average in the formula provided by Theorem~\ref{thc}
is now equal to $(\pi^2 d^{n+2})^{n-r} $ times 
\beq\label{bbe}
\mathbb E \Big[	{\bf 1}_{\{\Ind (\pi(x,\alpha))=i\}} 
\Big| \det{}_\R \Big(
\langle Y_{|K},\epsilon\rangle-
a\frac{\nabla\pi_\C(x)_{|K}(1+O(\frac1d))}{\|\pi_\C(x)\| \pi^{\frac32} d^{\frac{2n+3}2}}
\Big)\Big|  \ |\  \ s(x)=0, \nabla s(x)=\alpha\Big].
\eeq
Recall that $ \pi(x,\alpha)$ defined by~(\ref{pix}).
Besides, by Corollary~\ref{coco},
$$\left(\langle Y_{|K},\epsilon\rangle \ |\ s(x)=0, \nabla s(x)=\alpha\right)\sim N\left(O\Big(\frac{\|(\beta,a,b)\|}{d^{\frac32}}\Big),
\Sigma_{GOE}^{n-r}\big(1+O(\frac1d)\big)\right),$$
where the constants are independent of $\alpha$ and $\epsilon$, and where $\Sigma_{GOE}^{n-r}$ denotes the covariance matrix $\Sigma_{GOE}$ defined by~(\ref{goe}) in dimension $n-r$. When $d$ grows to infinity, the average~(\ref{bbe}) is uniformly bounded above by an integrable map, since the pole generated by $\|\pi_\C(x)\|$ is integrable. Consequently, the dominated convergence theorem implies that
\beqe
 	\frac{1}{d^n\vol (U)}\mathbb E \#( \crit_i^p\cap U) \underset{d\to \infty}{\to}&& 
 	\frac{2^{r^2 -nr -2r+2}}{\pi^{(r-1)(n-r+1)}}
 	\vol (\grass_\C(n-r,n-1)\\&&
\vol (\grass_\C(r-1,r) )\\
&&
\mathbb E \left(|\det \beta |^{2(n-r+2)}\right)
\mathbb E \left({\bf 1}_{\{\Ind A=i\}} 
|\det A|^2\right),
\eeqe
where
 $A \in M_{n-r}(\C)$ has covariance $\Sigma^{n-r}_{GOE}$ and where we used the determinant equality~(\ref{dety}).
Note that we passed from the real determinant $\det_\R$ to the complex one for the random complex matrix $A$. 
Since its index is always $n-r$, all the averages divided by $d^n$ for
$i\neq n-r$ converge to $0$.
The computations of the expectations and volume are given in~\cite[Remark 3.1.1, proof of Theorem 3.5.1]{gayet2015expected}:
\beqe
\vol (\grass_\C(n-r,n-1) &=& \pi^{(n-r)(r-1)}
\frac{\prod_{j=1}^{r-1}\Gamma(j)}
{\prod_{j=n-r+1}^{n-1}\Gamma(j)}\\
\vol (\grass_\C(r-1,r) ) &=& \pi^{r-1}\frac{1}{\Gamma(r)}\\
\mathbb E |\det \beta |^{2(n-r+2)}& =&
2^{(r-1)(n-r+2)}\frac{\prod_{j=n-r+3}^{n+1}\Gamma(j)}
{\prod_{j=1}^{r-1}\Gamma(j)}\\
\mathbb E \left(|\det Y|^2 \right)& =& 2^{n-r}(n-r+1)!
\eeqe
The powers of 2 in the latter equalities come from different choices of the measures, more precisely our choice of the half in the exponentials. 
Hence,
\beqe
\frac{1}{d^n\vol (U)}\mathbb E \#( \crit_{n-r}^p\cap U) \underset{d\to \infty}{\to} n!{n-1 \choose r-1}.
\eeqe
\epr	
We give now a sketch proof of the affine case.	
\bpr[ of Theorem~\ref{BF0}]
Let $f=(f_1, \cdots, f_r)\in \C^r$ be the random Bargmann-Fock field. 
 For any $R>0$, 
let $p_R= p(\frac{\cdot}{R}),$
so that the associated complexification $ \pi_{R,\C}$ of $dp_R(x)$ satisfies $\pi_{R,\C}(x)= \frac{1}R \pi_\C(\frac{x}R).$
Note that $p_R$ is a Morse function on $R U$. 
By Corollary~\ref{coaff}, the field $X$ defined by~\ref{x} is non-degenerate, so that 
we can apply Theorem~\ref{thc}  on the open set $RU$. 
By the independance of the triplet $(f,\nabla_0 f,\nabla_0^2 f)$, 
the conditional expectation in Theorem~\ref{thc} equals
$$\mathbb E \Big[	{\bf 1}_{\{\Ind (\pi(x,\alpha))=i\}} 
\Big| \det \Big(
\langle \nabla^2_0 f(x)_{|K},\epsilon\rangle - 
\frac1R\frac{d\pi_\C(\frac{x}R)_{|K}}{\|\pi_\C(\frac{x}R)\|}\langle \alpha(h(x)),\epsilon(\mu)\rangle
\Big)\Big|  \Big].$$
Recall that $ \pi(x,\alpha)$ defined by~(\ref{pix}).
We make the change of variables
$(\beta, a, b)= \frac{1}{\sqrt \pi}\alpha_{|K^\perp}$ (as~(\ref{bab}))
and $Y=\frac{1}\pi \nabla_0^2 f_{|K}$,
and then the change of variables $y=x/R$. 
By Corollary~\ref{coaff}, we obtain
\beqe	\frac{1}{R^{2n}\vol (U)} \mathbb E \# (R U\cap \crit_i^p)&\underset{R\to +\infty}{\to}& 
\frac{2^{r^2 -nr -2r+2}}{\pi^{(r-1)(n-r+1)}}
	\vol (\grass_\C(n-r,n-1)\\
&&\vol (\grass_\C(r-1,r) )\\
&&\mathbb E \left(|\det \beta |^{2(n-r+2)}\right)
\mathbb E \left({\bf 1}_{\{\Ind A=i\}} 
|\det A|^2\right).
\eeqe
we conclude as in the projective case. 
\epr
\subsection{The boundary case}
In this paragraph, we apply Proposition~\ref{kac} to estimate
the mean number of critical points of the restriction of $p$ on the boundary of $Z_s$ inside $\partial U$,
where $U\subset M$ is an open set with smooth boundary and $M$ is complex. We begin by a description of 
the mixed complex geometry on the boundary of $U$. 

\paragraph{Complex geometry on the boundary.}
In the sequel, for any $x\in M$ and any real subspace $L\subset T_x M$, we denote by $L_\C$ the largest complex subspace in $L$. 
Let $U\subset M$ be a codimension 0 open set with smooth boundary $\partial U$. 
 \begin{definition}\label{mobo} Let $Z$ be a smooth manifold of dimension $m$, with $C^2$ boundary, and $p : Z\to \R$ a smooth function. Then, $p$ is said to be \emph{Morse} if
	there is no critical point on $\partial Z$, 
	if $p$ is Morse and if $p_{|\partial Z}$ is Morse a well.
\end{definition} 
Let $p: M\to \R$ be a Morse function, such that $p_{|U}$ 
is Morse in the sens of Definition~\ref{mobo}.
Let $H=\ker p_\partial\subset T\partial U. $
For any $x\in \partial U$ which is not a critical point of $p_\partial$, $\dim H=2n-2$. Moreover,
either $H=H_\C$ and in this case $\dim_\C H = n-1$, 
or $\dim_\C H_\C = n-2.$
The first situation is non-generic, but our result holds in this case as well.
We define 
$$ F_\partial = E_{|\partial U}\oplus \mathcal L^\C(TM,E)_{|T(\partial U)},$$
$$ W=\{(x,0,\alpha) \in F_{\partial}, \alpha \text{ onto and }\ker \alpha\subset \ker dp_\partial (x)\},$$
$W_x$ its fiber over $x\in M$, and  $$X(x) = (x,s(x), \nabla_\partial s(x))\in F_{\partial},$$
where $\nabla_\partial=\nabla_{|T\partial U}$ denotes the restriction of the Chern connection $\nabla$ on $E$ to the tangent space of the boundary of $U$.
For any $x\in \partial U$ and any $\alpha\in \mathcal L^\C_{onto} (T_xM,E_x)$, 
let $$K=K(x,\alpha)=\ker \alpha \cap T\partial U.$$
If $K\neq \ker \alpha$, then $\dim_\R K = 2n-2r-1$ and $\dim_\C K_\C = n-r-1.$
Assume now that $K\subset H$. Then, 
$K_\C \subset H_\C$. 
Let $g\in H$ be a (one of the two) unit vector such that $$K=K_\C \obot \R g.$$ 
Note that $\ker\alpha = K\oplus \R Jg$,
where $J$ denotes the complex structure $J: TM\to TM$. 
Now, $\dim_\R K^\perp = 2r$, where $\perp$ stands for the metric on $T\partial U$.
Moreover, $$ \dim_\R (K^\perp \cap H) = 2r-1$$
so that $ \dim_\C (K^\perp\cap H)_\C = r-1.$
Let $v\in H$ a unit vector such that 
$$K^\perp \cap H=(K^\perp\cap H)_\C \obot \R v.$$
Note that $ K^\perp  =  (K^\perp\cap H)_\C +  \R v+ \R h$,
where $h\in H^\perp\setminus\{0\}.$ 
Now,
let 
$$
\mu = \ker dp_\partial (x) (\alpha_{|K^\perp})^{-1}= \alpha(K^\perp\cap H)\subset E_x.
$$
and 
$
\mu_\C =  \alpha((K^\perp\cap H)_\C).
$
Finally, let $\epsilon$ be a unit vector in $\mu^\perp\subset E_x.$

\begin{lemma}\label{gur}Under the setting above, for any $x\in \partial U$,
	the real dimension of $W_x$ equals $2nr-2n+2r+1.$
\end{lemma}
\bpr
For any $(x,0,\alpha)\in W$,
$$ T_{(0,\alpha)}W_x = 
\{(0,\beta)\in E_x\times 
\mathcal L^\C (T_x M, E_x), \ 
\ker dp_\partial(x) (\alpha_{|K^\perp})^{-1}\beta_{|K} =0\}.$$
Since $\beta_{|K_\C}$ 
is a complex linear map, its image in $E_x$ is a complex subspace, so that 
$(0,\beta)\in W_x$ if and only if 
$$ \beta(K_\C) \subset \mu_\C \text{ and }
\langle \beta_{|\R g},\epsilon\rangle = 0.$$ 
Since $\dim_\C \alpha((K^\perp\cap H)_\C)= r-1$, 
the real dimension of $W_x$ equals
$$ \dim_\R W_x = 2r^2 +2(n-r-1)(r-1)+(2r-1), $$
where the first term equals $\dim_\R \mathcal L^\C (\ker^\perp \alpha, E_x)$ (here $\perp$ stands for $T_xM$), the second equals $\dim_\R \mathcal L^\C (K_\C, \mu_\C)$ and the third
equals $\dim_\R \{\beta \in \mathcal L^\C(g^\C,E_x), \langle \beta,\epsilon\rangle = 0\},$
where $g^\C=\R g+\R Jg$ denotes the complex line generated by $g$.
\epr
\paragraph{The projective case. }
We first specialize this setting to the projective setting of Theorem~\ref{main}.
Recall that the natural scale for the random sections of degree $d$ is $d^{-\frac12}$. Since the dimension of $\partial U$ is $2n-1$, we can guess that the average number of critical points of $p_{|\partial U\cap Z_s}$  should be bounded by $O(d^{\frac{2n-1}2}).$
\begin{proposition}\label{b0} 
	Let $n\geq 2$ and $1\leq r\leq n-1$  be integers, $M$ be a compact smooth K\"ahler manifold and $(L,h)$ be an ample complex line bundle over $M$, with curvature form $\omega$, $(E,h_E)$ be a holomorphic rank $r$ vector bundle and
	let		$U\subset M$ be a 
	$0$-codimension submanifold with smooth boundary. Then, for any Borel subset $A\subset \partial U$, 
	\begin{eqnarray*}
		\forall 0\leq i\leq 2n-2r-1, \ \frac{1}{d^{n-\frac12}}\mathbb E \#  \crit_i^{p_{|\partial U}} 	= O_{d\to\infty}(1).
	\end{eqnarray*}	
	Here the probability measure is the Gaussian one given by~(\ref{mesure}). 
\end{proposition}
\bpr
Since we only need a bound for the averages and not their exact asymptotics, 
we apply Theorem~\ref{kac} which is easier to handle with than Corollary~\ref{coro1}.
By Theorem~\ref{kac}, we have that for any	Borel subset $A\subset \partial U$, 
\beq\label{cocor}	\mathbb E \left[\# \left(A\cap \crit_i^{p_\partial}\right)\right] &=& \int_{x\in A}
\int_{\substack{\alpha\in \mathcal L_{onto}^\C(T_xM,E_x)_{|T_x\partial U}\\ \ker \alpha\subset \ker dp_\partial(x)}} 
\left|\det \alpha_{|\ker^\perp \alpha}	\right| \\ \nonumber
&&\mathbb E \Big[	{\bf 1}_{\{\Ind (\pi(x,\alpha))=i\}} 
\Big| \det \Big(
\langle \nabla_\partial^2 s(x)_{|\ker \alpha},\epsilon(x,\alpha)\rangle \\ \nonumber
&& - \frac{\nabla^{2}p_\partial(x)_{|\ker \alpha}}{\|dp_\partial(x)\|}\langle \alpha(h(x)) ,\epsilon(x,\alpha)\rangle
\Big)\Big|
 \left. \middle| \  s(x)=0, \nabla_\partial s(x)=\alpha\right.\Big]\\ \nonumber &&
\rho_{X(x)}(0,\alpha)
d\vol(\alpha)  d\vol(x),
\eeq
where $\rho_{X(x)}$ is the Gaussian density of $X(x)$ and $\perp $ refers to the orthogonality in $T\partial U$.  
As in the proof for  projective manifold case, in equation~(\ref{cocor})
we perform the change of variables $\beta = d^{\frac{n+1}2}\alpha$ and 
$Y=d^{\frac{n+2}2} \nabla^2 s$. Then, thanks to Lemma~\ref{gur} which provides the power of $d$ which pops up from $\vol (\alpha)$, the average equals $d^{n-\frac{1}2}$ times a multiple integral which converges to a convergent integral independent of $d$. 
\epr
\paragraph{The affine setting. }
For the Bargmann-Fock field, we have the similar proposition:
\begin{proposition}\label{BF2}Let $1\leq r\leq n $ be integers, $ f : \C^n\to \C^r$ be $r$ independent copies of the Bargmann-Fock field~(\ref{bafo}), $U\subset \C^n$ be an open subset with smooth boundary, and $p : \bar U$ be a smooth Morse function, in the sense of Definition~\ref{momomo}. Then,
	\beq \forall 0\leq i\leq 2n-2r-1,\ 
	\ \frac{1}{R^{2n-1}}\mathbb E \#\crit_i^{p_{|\partial (RU)}} = O_{R\to +\infty}(1).
	\eeq
\end{proposition}
\bpr This is very similar to the projective setting.\epr
\subsection{Proof of the main theorems}

Theorem~\ref{main} is a simple consequence of Theorem~\ref{projective} and Proposition~\ref{b0}. Indeed, Morse inequalities for manifolds with boundary hold:  
 \begin{theorem}\label{morse}(see~\cite[Theorem A]{laudenbach}). Under the setting of Definition~\ref{mobo}, assume hat $\bar Z$ is compact. 
 	 For any $i\in \{0,\cdots, m-1\}, $ let  $N_i$ be the number of boundary critical points of $p_{|\partial Z}$ of index $i$, such that $p$ increases in the direction of $Z$.
 	Then,
 	\begin{itemize}
\item (weak Morse inequalities)
$ \forall 0\leq i\leq m, \ 
b_i(Z) \leq \# \crit_i^p +N_i.$
\item (strong Morse inequalities)
$\forall 0\leq i\leq m, $
$$\sum_{k=0}^i (-1)^{i-k}b_k(Z)\geq 
\sum_{k=0}^i (-1)^{i-k} (\# \crit_k^p +N_k).$$
\end{itemize}
 \end{theorem}
 We will apply these Morse inequalities to the random nodal sets $Z_s \cap U$.

 \bpr[ of Theorem~\ref{main}] 
 By~\cite[Lemma 2.8]{GWdet}, almost surely the restriction $p_{|Z_s}$ is Morse in the latter sense. The proof of this lemma extends to $p_{|Z_s\cap U}$, so that we can apply Theorem~\ref{morse} to $Z_s\cap U$, for almost all $s$. Hence,
\beq\label{momomo}
 \forall 0\leq i\leq 2n-2r, \ 
\mathbb E b_i(Z_s\cap U) \leq \mathbb E (\# \crit_i^p) + \mathbb E (\# \crit_i^{p_{|\partial U}}).
\eeq
By (\ref{momomo}), Theorem~\ref{projective} and Proposition~\ref{b0}, we obtain
$$ \forall 0\leq i\leq 2n-2r\setminus\{n-r\}, \ 
\mathbb E b_i(Z_s\cap U)=o(d^n)$$
and 
$\mathbb E b_{n-r}(Z_s\cap U)\leq  \mathbb E (\#\crit_{n-r}^p)  +o(d^n).$
On the other hand, the two assertions of Theorem~\ref{morse} and Proposition~\ref{b0} imply that  
$$\mathbb E b_{n-r}(Z_s\cap U)\geq  \mathbb E (\#\crit_{n-r}^p)  -o(d^n),$$
so that by Theorem~\ref{projective},
$ \mathbb E b_{n-r} = d^n{n-1\choose r-1}\int_U\omega +o(d^n),$
which is the result.
\epr
\begin{lemma}\label{wig}(\cite[Lemma 3.2]{wigman2021expected})
	Under the hypotheses of Theorem~\ref{BF}, there exists a map  $p: \C^n \to \R$ such that for almost all instance of the Bargmann-Fock field, $p_{|Z_s\cap \bar U}$ is Morse as well in the sense of Definition~\ref{mobo}.
\end{lemma}
\bpr This is proven (in a more general setting) in the proof of \cite[Lemma 3.2]{wigman2021expected} for a manifold without boundary. The argument extends immediatly to manifolds with $C^2$ boundary. \epr
\bpr[ of Theorem~\ref{BF}] The proof is similar to the one of Theorem~\ref{main}, using Theorem~\ref{BF0} and Proposition~\ref{BF2}.
\epr

\bibliographystyle{amsplain}
\bibliography{coursgeo.bib}

\noindent Univ. Grenoble Alpes, Institut Fourier \\
F-38000 Grenoble, France \\
CNRS UMR 5208  \\
CNRS, IF, F-38000 Grenoble, France

\end{document}